\documentclass[11pt]{article}
\usepackage[T1]{fontenc}
\usepackage{lmodern}
\usepackage{amsmath,amssymb,amsthm}
\usepackage{graphicx}
\usepackage{geometry}
\usepackage{microtype}
\usepackage{caption}
\usepackage{xcolor}
\usepackage[colorlinks=true,linkcolor=blue,citecolor=red,urlcolor=blue]{hyperref}
\hypersetup{pdftitle={Preservation and Loss of Concavity in the Vlasov--Poisson System},pdfauthor={Zhaopeng Wang}}
\numberwithin{equation}{section}
\newtheorem{theorem}{Theorem}[section]
\newtheorem{lemma}[theorem]{Lemma}
\newtheorem{remark}[theorem]{Remark}

\newcommand{\R}{\mathbb R}

\newcommand{\la}{\langle}
\newcommand{\ra}{\rangle}
\title{Preservation and Loss of Log-Concavity in the Repulsive Vlasov--Poisson System}
\author{
Xiaoxiao Sun\textsuperscript{1}
\qquad
Zhaopeng Wang\textsuperscript{2}\\[0.4em]
School of Mathematics and Statistics\\
Beijing Institute of Technology\\
Beijing 100081, China\\[0.3em]
\textsuperscript{1}\href{mailto:xsun@bit.edu.cn}{xsun@bit.edu.cn}
\qquad
\textsuperscript{2}\href{mailto:zpwang@bit.edu.cn}{zpwang@bit.edu.cn}
}
\date{}
\begin{document}
\maketitle

\begin{abstract}
Pr\'ekopa's theorem shows that, for jointly log-concave initial data, free transport preserves the  log-concavity of spatial densities.  This paper asks whether this property persists in the repulsive Vlasov--Poisson system. In dimensions two and higher, we construct smooth, strongly log-concave initial distributions whose spatial densities lose log-concavity at arbitrarily small positive times. For one-dimensional cold-ion data, we prove that the spatial density remains log-concave for all time.
\end{abstract}



\section{Introduction}

Log-concavity is a simple but useful property of functions. A positive function $g$ is called log-concave if $\log g$ is concave. Gaussian functions, for example, are log-concave.  A classical theorem of Pr\'ekopa \cite{Prekopa1973} says that integration with respect to some of the variables preserves log-concavity: if $g(x,v)$ is jointly log-concave, then
\[
 x\longmapsto\int_{\R^d}g(x,v)\,dv
\]
is also log-concave, provided that the integral is finite. This result is a direct consequence of the Pr\'ekopa--Leindler inequality, see  Brascamp and Lieb \cite{BrascampLieb1976} for a simple proof and further developments. For background on log-concavity and strong log-concavity, we refer to the survey of Saumard and Wellner \cite{SaumardWellner2014}.

This result is particularly natural in kinetic theory, where the spatial density is obtained by integrating the particle distribution over velocity. Consider first the free transport equation
\[
 \partial_t f+v\cdot\nabla_x f=0,
 \qquad f(0,x,v)=f_0(x,v).
\]
Its solution is
\[
 f(t,x,v)=f_0(x-tv,v).
\]
For each fixed $t\geq 0$, the map $(x,v)\mapsto(x-tv,v)$ is affine. Hence joint log-concavity is preserved, and Pr\'ekopa's theorem gives the log-concavity of the spatial
density for every $t\geq0$.

\begin{figure}[tbp]
 \centering
 \includegraphics[width=.91\textwidth]{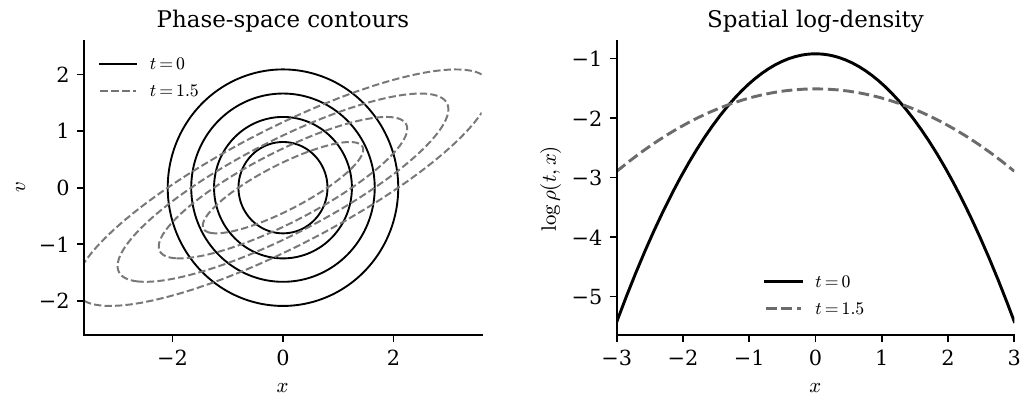}
 \caption{A Gaussian under free transport and its log-concave spatial marginal.}
\end{figure}

The reason is simple: under free transport, the particles do not interact. With a self-consistent field, each particle is subject to the force generated by the spatial density. The density determines the field, and the field changes the particle motion in return. Consequently, the characteristic flow is generally no longer affine. This raises the question of whether the spatial density remains log-concave?

The Vlasov--Poisson system is a classical model for such self-consistent motion. In the repulsive case, it takes the form
\begin{equation}\label{eq:VP}
\begin{cases}
 \partial_t f+v\cdot\nabla_x f+E\cdot\nabla_v f=0,\quad f(0,x,v)=f_0(x,v),\\
 \nabla_x\cdot E=\rho,\qquad
 \displaystyle\rho(t,x)=\int_{\R^d}f(t,x,v)\,dv.
\end{cases}
\end{equation}
Here $t\geq0$, $x,v\in\R^d$,  $f(t,x,v)\geq0$ is the particle distribution, and $\rho(t,x)$ is the spatial density. The electric field is the repulsive Coulomb field
\begin{equation}\label{eq:Coulomb}
 E(t,x)=\frac1{|\mathbb S^{d-1}|}
 \int_{\R^d}\frac{x-y}{|x-y|^d}\rho(t,y)\,dy.
\end{equation}
When $d=1$,  formula \eqref{eq:Coulomb} becomes
$$
E(t,x)=\int_{-\infty}^x\rho(t,y)\,dy-\frac12\int_{\R}\rho(t,y)\,dy.
$$

The classical theory of the Vlasov--Poisson system is by now well developed. Global classical solutions in one space dimension go back to Iordanskii \cite{Iordanskii1961}; the two-dimensional case was treated by Ukai and Okabe \cite{UkaiOkabe1978}. In three dimensions, Batt proved the existence of global symmetric solutions \cite{Batt1977}, while Bardos and Degond proved the existence of global smooth solutions for small initial data \cite{BardosDegond1985}. Global existence of classical solutions for general initial smooth data was subsequently established by Pfaffelmoser \cite{Pfaffelmoser1992}, Schaeffer \cite{Schaeffer1991}, and Lions and Perthame \cite{LionsPerthame1991}. For a broader account of the well-posedness theory for Vlasov-type equations, including the classical Vlasov–Poisson system and related ion models, we refer to  \cite{GriffinPickeringIacobelli2021b, Rein2007}. For ion models with massless electrons, Griffin-Pickering and Iacobelli also proved global strong well-posedness for two ionic Vlasov–Poisson systems in \(\mathbb R^3\)  \cite{GriffinPickeringIacobelli2021a}.

In this paper, we consider a different question: whether the log-concavity of the spatial density is preserved in the repulsive Vlasov--Poisson system? In dimensions $d\geq2$, the answer can be negative, even for very regular initial data.

\begin{theorem}\label{thm:high-d-intro}
Let $d\geq2$ and $\tau>0$. There exist $c>0$ and an initial distribution $f_0$ such that
\[
 0<f_0\in C^\infty(\R^{2d})\cap L^1(\R^{2d}),
 \qquad -D^2_{x,v}\log f_0\geq cI_{2d}.
\]
For the corresponding local classical solution of \eqref{eq:VP}, there exist $t_*\in(0,\tau)$ and $x_*\in\R^d$ such that
\[
 \lambda_{\max}\bigl(D_x^2\log\rho(t_*,x_*)\bigr)>0.
\]
\end{theorem}

The initial distribution may be chosen radially symmetric in $x$ and Gaussian in $v$. Since $\tau>0$ is arbitrary, loss of log-concavity can occur at arbitrarily small positive times. However, the initial datum may depend on $\tau$; thus, we do not claim that one fixed datum loses log-concavity  at arbitrarily small times.

The self-consistent field does not always destroy log-concavity. We next consider particles that all start at rest. Their initial distribution is
\begin{equation}\label{eq:cold-data}
 f_0(x,v)=\rho_0(x)\delta_0(v).
\end{equation}
This is a measure in velocity. In one dimension, it is related to the pressureless Euler--Poisson dynamics studied, for example, by Engelberg, Liu, and Tadmor \cite{EngelbergLiuTadmor2001} and by Carrillo and Galtung \cite{CarrilloGaltung2026}.

\begin{theorem}\label{thm:cold-intro}
Let $d=1$ and suppose that
\[
 0<\rho_0\in C^\infty(\R)\cap L^1(\R),
 \qquad (\log\rho_0)''\leq0.
\]
Then, for the initial data \eqref{eq:cold-data}, there exists a global measure-valued solution to \eqref{eq:VP} whose spatial density satisfies
\[
 \partial_x^2\log\rho(t,x)\leq0,
 \qquad t\geq0,\quad x\in\R.
\]
\end{theorem}

\begin{remark}
It is also natural to ask whether the loss and preservation properties of log-concavity established in this paper have counterparts for the attractive Vlasov--Poisson equation, obtained from \eqref{eq:VP} by replacing \(+E\cdot\nabla_v f\) with \(-E\cdot\nabla_v f\), while keeping the Poisson equation unchanged. Since the force has the opposite sign, this problem may require a different mechanism.
\end{remark}

The proofs rely two different features of \eqref{eq:VP}. To prove  Theorem~\ref{thm:high-d-intro}, we compute the second time derivative of the radial logarithmic curvature at $t=0$. A radial convex potential with a linear tail makes this derivative positive at a suitable point. We then add a small quadratic term to the potential, making the initial distribution becomes strongly log-concave. The stability estimate in Appendix~A ensures that the positive curvature remains at a fixed small positive time.

The proof of Theorem~\ref{thm:cold-intro} uses an explicit characteristic map. Let $X(t,a)$ denote the position at time $t$ of the particle starting at $a$. Then
\[
 X(t,a)=a+\frac{t^2}{2}E_0(a),
 \qquad
 E_0(a)=\int_{-\infty}^a\rho_0(s)\,ds-\frac12\int_{\R}\rho_0(s)\,ds.
\]
Since $\partial_aX(t,a)=1+\frac{t^2}{2}\rho_0(a)>0$, the particle order is preserved. Conservation of mass then yields 
\[
 \rho(t,X(t,a))=
 \frac{\rho_0(a)}{1+\frac{t^2}{2}\rho_0(a)}.
\]
A direct differentiation of this formula gives the desired log-concavity.

The two theorems concern different classes of solutions. Theorem~\ref{thm:high-d-intro} concerns smooth phase-space distributions in $d\geq2$. Theorem~\ref{thm:cold-intro} concerns a one-dimensional cold-ion flow. Thus, these results do not give a dimension-only dichotomy: we make no claim here for smooth data in $d=1$ or for cold-ion data in $d\geq2$. Together they show that a self-consistent field may destroy log-concavity, but need not do so.

Section~2 gives the construction leading to the loss of log-concavity and then proves Theorem~\ref{thm:high-d-intro}. Section~3 proves Theorem~\ref{thm:cold-intro}. The local existence and stability estimates used in Section~2 are established in Appendix~A.

\section{\texorpdfstring{$d\geq2$}{d >= 2}}

\subsection{The basic construction}

We study smooth solutions of \eqref{eq:VP}, with the Coulomb field defined by \eqref{eq:Coulomb}. The data used below are smooth and have rapid decay in both variables. The local existence and time regularity needed for the computations are stated in Lemma~\ref{lem:stability} and proved in Appendix~A.

Let
\[
 j(t,x)=\int_{\R^d}v f(t,x,v)\,dv,
 \qquad
 \Pi(t,x)=\int_{\R^d}v\otimes v f(t,x,v)\,dv.
\]
Integrating \eqref{eq:VP} against $1$ and $v$ gives
\begin{equation}\label{eq:moments}
 \partial_t\rho+\nabla_x\cdot j=0,
 \qquad
 \partial_tj+\nabla_x\cdot\Pi=\rho E.
\end{equation}
We choose a centered Gaussian in velocity:
\begin{equation}\label{eq:maxwell}
 f_0(x,v)=\rho_0(x)(2\pi\theta)^{-d/2}
 \exp\!\left(-\frac{|v|^2}{2\theta}\right),
 \qquad \theta>0.
\end{equation}
At $t=0$, $j_0=0$ and $\Pi_0=\theta\rho_0 I_d$. It follows from
\eqref{eq:moments} that
\[
 \partial_t\rho(0)=0,
 \qquad
 \partial_t^2\rho(0)=\theta\Delta\rho_0-\nabla_x\cdot(\rho_0E_0).
\]
Writing $u=\log\rho$ and $u_0=\log\rho_0$, we obtain
\begin{equation}\label{eq:curvature-source}
 \partial_tu(0)=0,
 \qquad
 \partial_t^2u(0)=
 \theta\bigl(\Delta u_0+|\nabla_xu_0|^2\bigr)
 -\rho_0-E_0\cdot\nabla_xu_0.
\end{equation}
This identity does not require radial symmetry. It suggests how to produce loss of log-concavity. We look for a point where the initial logarithmic curvature is zero but its second time derivative is positive. Taylor's formula then gives positive curvature at small positive times. We now construct data with this property.

Choose $b\in C_c^\infty([0,\infty))$ such that $b\geq0$, $b$ is a positive constant near $0$, and $b$ is not identically zero. Set
\begin{equation}\label{eq:potential}
 V(r)=\int_0^r(r-s)b(s)\,ds,
 \qquad r\geq0.
\end{equation}
Then $V''=b\geq0$ and $V'(r)=\int_0^r b(s)\,ds\geq0$. The radial and tangential eigenvalues of $D_x^2V(|x|)$ are $V''(r)$ and $V'(r)/r$.
Both are nonnegative. Since $b$ is constant near the origin, $V(|x|)$ is quadratic there and is smooth on all of $\R^d$. Since $b$ has compact support, there are $R_0>0$, $C_0\in\R$, and $\lambda>0$ such that
\[
 V(r)=\lambda r+C_0,\qquad r\geq R_0.
\]
Thus $V(|x|)$ is convex and has a linear tail.

Let $q(r)=e^{-V(r)}$ and set $\rho_0(x)=Aq(|x|)$, where $A>0$ will be chosen later. The initial field is radial, $E_0(x)=e(r)x/r$, and Gauss' law gives
\begin{equation}\label{eq:eradial}
 e(r)=\frac{A}{r^{d-1}}\int_0^r s^{d-1}q(s)\,ds,
 \qquad
 e'(r)=Aq(r)-\frac{d-1}{r}e(r).
\end{equation}
For $r\geq R_0$, we have
\[
 u_0(r)=\log A-\lambda r-C_0,
 \qquad \Delta u_0=-\frac{(d-1)\lambda}{r},
 \qquad |\nabla_xu_0|^2=\lambda^2.
\]
Therefore \eqref{eq:curvature-source} becomes
\[
 \partial_t^2u(0,r)
 =\theta\left(\lambda^2-\frac{(d-1)\lambda}{r}\right)
 -Aq(r)+\lambda e(r).
\]
Differentiating the second identity in \eqref{eq:eradial} and using $q'(r)=-\lambda q(r)$ gives
\[
 e''(r)=-A\left(\lambda+\frac{d-1}{r}\right)q(r)
 +\frac{d(d-1)}{r^2}e(r).
\]
Using $q''(r)=\lambda^2q(r)$ and the formula for $e''$, a second derivative gives
\begin{equation}\label{eq:radial-curvature}
 \partial_r^2\partial_t^2u(0,r)
 =\frac{\lambda d(d-1)}{r^2}e(r)
 -A\left(2\lambda^2+\frac{(d-1)\lambda}{r}\right)q(r)
 -\frac{2\lambda(d-1)\theta}{r^3}.
\end{equation}
The first term has an algebraic tail, while the second decays exponentially. The last term is the contribution of the velocity spread.
The initial distribution is invariant under simultaneous orthogonal transformations of $x$ and $v$. Uniqueness and invariance of \eqref{eq:VP} give $f(t,Qx,Qv)=f(t,x,v)$ for $Q\in O(d)$.
Integrating in $v$ gives $\rho(t,Qx)=\rho(t,x)$, so $\rho$ remains radial. Hence $u_{rr}$ is the radial eigenvalue of $D_x^2\log\rho$.

We now make this comparison precise. Put
\[
 h(r)=r^{1-d}\int_0^r s^{d-1}q(s)\,ds,
 \qquad
 M_q=\int_0^\infty s^{d-1}q(s)\,ds>0.
\]
Then $e(r)=Ah(r)$ and
\[
 h(r)=\frac{M_q+o(1)}{r^{d-1}},
 \qquad q(r)=e^{-\lambda r-C_0}
 \quad\text{as }r\to\infty.
\]
Define
\[
 F(R):=\frac{\lambda d(d-1)}{R^2}h(R)
 -\left(2\lambda^2+\frac{(d-1)\lambda}{R}\right)q(R).
\]
Since $h(R)=M_qR^{1-d}(1+o(1))$ and $q(R)=e^{-\lambda R-C_0}$,
\[
 F(R)=\lambda d(d-1)M_qR^{-(d+1)}(1+o(1))
 -\left(2\lambda^2+\frac{(d-1)\lambda}{R}\right)e^{-\lambda R-C_0}.
\]
The algebraic term dominates the exponential term. Thus $F(R)>0$ for all sufficiently large $R>R_0$. Fix such an $R$. We then choose
\[
 A>\frac{2\lambda(d-1)\theta}{R^3F(R)},
\]
so that
\begin{equation}\label{eq:positive-B}
 B:=AF(R)-\frac{2\lambda(d-1)\theta}{R^3}>0.
\end{equation}
By \eqref{eq:radial-curvature}, $\partial_r^2\partial_t^2u(0,R)=B$. At the chosen point, $u_0''(R)=0$, and \eqref{eq:curvature-source} gives $u_t(0,r)=0$ for every $r$. Hence Taylor's formula and
\eqref{eq:radial-curvature} show that
\begin{equation}\label{eq:loss-flat}
 \partial_r^2\log\rho(t,R)
 =\frac{B}{2}t^2+o(t^2)>0
 \qquad\text{for }0<t<T_1,
\end{equation}
for some $T_1>0$. Since $B>0$, the radial eigenvalue of $D_x^2\log\rho(t,Re_1)$ is positive for $0<t<T_1$. This gives a smooth, jointly log-concave initial distribution whose spatial density loses log-concavity.

\subsection{Proof of Theorem 1.1}

The construction gives a jointly log-concave initial distribution. We now make it strongly log-concave without losing the positive curvature at positive time. The only additional ingredient is continuity of the density in the local $C^2$ norm. We state the required stability result first.

For an integer $s\geq0$ and $m>0$, define
\begin{equation}\label{eq:weighted-norm}
 \|g\|_{H^s_m}^2
 :=\sum_{|\alpha|+|\beta|\leq s}
 \left\|\la v\ra^m
 \partial_x^\alpha\partial_v^\beta g
 \right\|_{L^2_{x,v}}^2,
 \qquad \la v\ra=(1+|v|^2)^{1/2}.
\end{equation}

\begin{lemma}\label{lem:stability}
Let $d\geq2$, let $s>d+4$ be an integer, and let $m>d/2+2$. We use the norm defined in \eqref{eq:weighted-norm}. Suppose that $\{g_{\rm in}^\varepsilon\}_{0\leq\varepsilon\leq1}$ is a family of nonnegative smooth data satisfying
\[
 \sup_{0\leq\varepsilon\leq1}
 \left(\|g_{\rm in}^\varepsilon\|_{H^s_m}
 +\|g_{\rm in}^\varepsilon\|_{W^{1,1}(\R^{2d})}\right)<\infty.
\]
Then the corresponding solutions to \eqref{eq:VP} exist on a common interval $[0,T]$ and satisfy
\[
 g^\varepsilon\in
 C([0,T];L^1\cap H^{s-1}_m)
 \cap L^\infty(0,T;H^s_m\cap W^{1,1}(\R^{2d})),
\]
with bounds independent of $\varepsilon$. Moreover, for any $\varepsilon,\eta\in[0,1]$,
\[
\begin{aligned}
 \sup_{0\leq t\leq T}
 \bigl(&\|g^\varepsilon(t)-g^\eta(t)\|_{L^1}
       +\|g^\varepsilon(t)-g^\eta(t)\|_{H^{s-1}_m}\bigr)\\
 &\leq C_T\bigl(
 \|g_{\rm in}^\varepsilon-g_{\rm in}^\eta\|_{L^1}
 +\|g_{\rm in}^\varepsilon-g_{\rm in}^\eta\|_{H^{s-1}_m}\bigr).
\end{aligned}
\]
In particular, if
\[
 g_{\rm in}^\varepsilon\longrightarrow g_{\rm in}^0
 \quad\text{in }H^{s-1}_m\cap L^1,
\]
then
\[
 \rho^\varepsilon\longrightarrow\rho^0
 \quad\text{in }C([0,T];C^2_{\rm loc}(\R^d)).
\]
\end{lemma}
\begin{remark}
For the smooth data used in this paper, the solution has enough time regularity for the second-order expansion in \eqref{eq:loss-flat}.
The proof is given in Appendix~A.
\end{remark}

\begin{proof}[Proof of Theorem~\ref{thm:high-d-intro}]
Fix $d\geq2$ and $\tau>0$. Choose $V$, $R$, and $A$ as in \eqref{eq:potential} and \eqref{eq:positive-B}. For $0\leq\varepsilon\leq1$, set
\begin{equation}\label{eq:strong-data}
 f_{0,\varepsilon}(x,v)
 =A\exp\!\bigl(-V(|x|)-\varepsilon|x|^2\bigr)
 (2\pi\theta)^{-d/2}
 \exp\!\left(-\frac{|v|^2}{2\theta}\right).
\end{equation}
The family has uniform bounds in the spaces of Lemma~\ref{lem:stability}. Indeed, outside a compact set $V(r)=\lambda r+C_0$, and the velocity factor is Gaussian. Every fixed derivative of the spatial factor is bounded, uniformly for $0\leq\varepsilon\leq1$, by a polynomial times $e^{-V(r)}$. This gives the required uniform weighted Sobolev and $W^{1,1}$ bounds. Moreover,
\[
 f_{0,\varepsilon}\longrightarrow f_{0,0}
 \quad\text{in }H^{s-1}_m\cap L^1
 \quad\text{as }\varepsilon\to0.
\]
This follows from dominated convergence: outside a compact set, $V(r)=\lambda r+C_0$, and every fixed derivative of $e^{-\varepsilon|x|^2}$ is bounded uniformly for $0\leq\varepsilon\leq1$ by a polynomial times the exponential spatial tail. Lemma~\ref{lem:stability} gives a common existence time $T>0$ and $\rho^\varepsilon\to\rho^0$ in $C([0,T];C^2_{\rm loc})$.

By \eqref{eq:loss-flat}, choose
\[
 0<t_*<\min\{\tau,T,T_1\}.
\]
Then
\[
 \gamma:=\partial_r^2\log\rho^0(t_*,R)>0.
\]
The unperturbed solution is radial by the observation in Section~2.1. Hence
\[
 \partial_{11}\log\rho^0(t_*,Re_1)=\gamma.
\]
The unperturbed density is positive at $(t_*,Re_1)$. Since $\rho^\varepsilon\to\rho^0$ locally uniformly, the perturbed densities have a common positive lower bound in a neighborhood of this point for all sufficiently small $\varepsilon$. The $C^2$ convergence implies
\[
 \partial_{11}\log\rho^\varepsilon(t_*,Re_1)
 \longrightarrow
 \partial_{11}\log\rho^0(t_*,Re_1)=\gamma.
\]
For sufficiently small $\varepsilon>0$, the left side is positive. On the other hand, convexity of $V(|x|)$ gives
\[
 -D_{x,v}^2\log f_{0,\varepsilon}
 =\begin{pmatrix}
 D_x^2V(|x|)+2\varepsilon I_d&0\\
 0&\theta^{-1}I_d
 \end{pmatrix}
 \geq\min\{2\varepsilon,\theta^{-1}\}I_{2d}.
\]
Choose such an $\varepsilon$ and set $f_0=f_{0,\varepsilon}$, $x_*=Re_1$, and $c=\min\{2\varepsilon,\theta^{-1}\}$. Then $f_0$ satisfies the assumptions of Theorem~\ref{thm:high-d-intro}, while its spatial density has positive logarithmic curvature at $(t_*,x_*)$.
\end{proof}

\section{\texorpdfstring{$d=1$}{d = 1}}

We now construct the measure solution of \eqref{eq:VP} with initial data \eqref{eq:cold-data}. The function $\rho_0$ satisfies the assumptions of Theorem~\ref{thm:cold-intro}. Let
\[
 M=\int_{\R}\rho_0(a)\,da.
\]
By the one-dimensional form of the Coulomb field in \eqref{eq:Coulomb},
\begin{equation}\label{eq:one-d-field}
 E_0(a)=\int_{-\infty}^a\rho_0(s)\,ds-\frac M2.
\end{equation}
Thus $E_0'=\rho_0$ and $|E_0|\leq M/2$. Define the particle paths
\[
 X(t,a)=a+\frac{t^2}{2}E_0(a),
 \qquad U(t,a)=tE_0(a).
\]
Their derivative in $a$ is
\[
 \partial_aX(t,a)=1+\frac{t^2}{2}\rho_0(a)>0.
\]
Since $E_0$ is bounded, $X(t,a)-a$ is bounded for each fixed $t$. Therefore $X(t,\cdot)$ is an increasing diffeomorphism of $\R$ onto itself for every $t\geq0$.

Push the measure $\rho_0(a)\,da$ forward by the map $a\mapsto(X(t,a),U(t,a))$. The resulting phase-space measure is
\[
 f(t,dx,dv)=\int_{\R}\rho_0(a)
 \delta_{X(t,a)}(dx)\delta_{U(t,a)}(dv)\,da.
\]
At $t=0$, this measure is exactly the initial datum \eqref{eq:cold-data}. Its spatial marginal is the pushforward by $X(t,\cdot)$. Hence
\begin{equation}\label{eq:mass-identity}
 \rho(t,X(t,a))\,\partial_aX(t,a)=\rho_0(a),
 \qquad
 \rho(t,X(t,a))=
 \frac{\rho_0(a)}{1+\frac{t^2}{2}\rho_0(a)}.
\end{equation}

We check that this measure solves \eqref{eq:VP}. Since $X(t,\cdot)$ is increasing, the mass to the left of $X(t,a)$ equals the initial mass to the left of $a$. Therefore, with the same normalization as in \eqref{eq:one-d-field},
\[
 E(t,X(t,a))=E_0(a).
\]
It follows that $\partial_tX=U$ and $\partial_tU=E(t,X)$. For a test function $\varphi\in C_c^\infty([0,\infty)\times\R^2)$, differentiation gives
\[
 \frac{d}{dt}\int_{\R}\rho_0(a)
 \varphi(t,X(t,a),U(t,a))\,da
 =\int_{\R}\rho_0(a)
 \bigl(\partial_t\varphi+U\partial_x\varphi
 +E(t,X)\partial_v\varphi\bigr)(t,X,U)\,da.
\]
Integrating in time and using the compact support of $\varphi$ gives
\[
 \int_0^\infty\!\int_{\R}\rho_0(a)
 \bigl(\partial_t\varphi+U\partial_x\varphi
 +E(t,X)\partial_v\varphi\bigr)(t,X,U)\,da\,dt
 +\int_{\R}\rho_0(a)\varphi(0,a,0)\,da=0.
\]
This is the weak formulation of the transport equation in \eqref{eq:VP}, including the initial trace. Thus the construction is self-consistent and gives a global measure solution. The spatial density is smooth and positive because $X(t,\cdot)$ is a smooth diffeomorphism and $\rho_0$ is smooth and positive.

\begin{proof}[Proof of Theorem~\ref{thm:cold-intro}]
Fix $t\geq0$ and write
\[
 q(a)=\rho_0(a),\qquad c=\frac{t^2}{2},\qquad
 P(a)=1+cq(a).
\]
By \eqref{eq:mass-identity},
\[
 \log\rho(t,X(t,a))=\log q(a)-\log P(a).
\]
Since $\partial_aX=P$, differentiation with respect to $x$ at $x=X(t,a)$ is given by $\partial_x=P^{-1}\partial_a$. Thus
\[
 \partial_x\log\rho(t,X(t,a))
 =\frac{q'(a)}{q(a)P(a)^2}.
\]
Differentiating once more, we obtain
\[
 \partial_x^2\log\rho(t,X(t,a))
 =\frac{(\log q)''(a)}{P(a)^3}
 -\frac{2c(q'(a))^2}{q(a)P(a)^4}.
\]
Both terms are nonpositive. Hence
\[
 \partial_x^2\log\rho(t,X(t,a))\leq0
 \qquad(t\geq0,\ a\in\R).
\]
Since $X(t,\cdot)$ maps $\R$ onto $\R$, this proves the desired inequality at every $x\in\R$.
\end{proof}

\appendix
\section{Local existence and stability}

We prove Lemma~\ref{lem:stability}. The local transport construction for Vlasov equations is classical; see Horst \cite{Horst1981}. We give the estimates needed for the stronger stability statement used in Section~2. Throughout this appendix $d\geq2$, $E[g]$ denotes the Coulomb field generated by
$\rho_g=\int g\,dv$, and $w_m(v)=\la v\ra^m$.

\subsection{Basic estimates}

For $m>d/2$, Cauchy--Schwarz in $v$ gives, for $0\leq k\leq s$,
\begin{equation}\label{eq:A-density}
 \|\rho_g\|_{H^k_x}\leq C_{d,m}\|g\|_{H^k_m}.
\end{equation}
Indeed, for each spatial derivative,
\[
 |\partial_x^\alpha\rho_g(x)|^2
 \leq\left(\int_{\R^d}\la v\ra^{-2m}\,dv\right)
 \left(\int_{\R^d}\la v\ra^{2m}
 |\partial_x^\alpha g(x,v)|^2\,dv\right),
\]
and integration in $x$ gives the estimate.

The Coulomb kernel satisfies $|K_d(z)|\leq C_d|z|^{1-d}$. Split its convolution into $|x-y|\leq1$ and $|x-y|>1$. The first part is bounded by $C\|\rho_g\|_\infty$, because $|z|^{1-d}$ is locally integrable. The second part is bounded by $C\|\rho_g\|_1$, because $|z|^{1-d}\leq1$ when $|z|>1$. Thus
\begin{equation}\label{eq:A-field-infty}
 \|E[g]\|_{L^\infty_x}
 \leq C_d\bigl(\|\rho_g\|_{L^1_x}
 +\|\rho_g\|_{L^\infty_x}\bigr).
\end{equation}
This estimate also holds for signed densities. In particular, it applies to differences of two solutions.

Each entry of $\nabla_xE[g]$ is a zero-order Fourier multiplier of $\rho_g$. The Calder\'on--Zygmund estimate therefore gives
\begin{equation}\label{eq:A-field-gradient}
 \|\nabla_xE[g]\|_{H^k_x}
 \leq C_{d,k}\|\rho_g\|_{H^k_x},
 \qquad k\geq0.
\end{equation}
Combining \eqref{eq:A-density}, \eqref{eq:A-field-infty}, and \eqref{eq:A-field-gradient}, we obtain
\[
 \|E[g]\|_\infty+\|\nabla E[g]\|_{H^{s-1}_x}
 \leq C\bigl(\|g\|_1+\|g\|_{H^s_m}\bigr).
\]
We do not use an $L^2$ bound for $E[g]$. Such a bound need not hold when the total charge is nonzero in dimension two.

\par\medskip
We next derive the weighted energy estimate. We first record the product estimate used for the commutators. Choose $\ell=\lfloor d/2\rfloor+1$. Then $\ell>d/2$ and
$s>2\ell+2$, since $s>d+4$. For smooth $E=E(x)$ and $g=g(x,v)$, Leibniz' rule and Sobolev embedding give
\begin{equation}\label{eq:A-product}
 \sum_{\substack{|\alpha|+|\beta|\leq s\\0<\gamma\leq\alpha}}
 \left\|w_m(\partial_x^\gamma E)
 \partial_x^{\alpha-\gamma}\partial_v^{\beta+1}g\right\|_2
 \leq C\|\nabla E\|_{H^{s-1}_x}\|g\|_{H^s_m}.
\end{equation}
Here and below, a sum over $\gamma$ is a sum over multi-indices. For completeness, split the terms into two cases. If
$|\gamma|\leq s-\ell$, then
\[
 \|\partial_x^\gamma E\|_\infty
 \leq C\|\partial_x^\gamma E\|_{H^\ell_x}
 \leq C\|\nabla E\|_{H^{s-1}_x}.
\]
The last inequality follows because $|\gamma|\geq1$ and $|\gamma|+\ell-1\leq s-1$. We put this factor in $L^\infty_x$ and the differentiated $g$ in
$L^2_{x,v}$. If $|\gamma|>s-\ell$, we use $\|\partial_x^\gamma E\|_2\leq\|\nabla E\|_{H^{s-1}_x}$ and put this factor in $L^2_x$. The remaining factor is put in
$L^\infty_x(L^2_v)$, using $H^\ell_x(L^2_v)\hookrightarrow L^\infty_x(L^2_v)$. In this case,
\[
 |\alpha-\gamma|+|\beta|+1+\ell
 \leq s-|\gamma|+1+\ell\leq 2\ell\leq s.
\]
Thus the required $\ell$ additional $x$-derivatives of $g$ are available.
This proves \eqref{eq:A-product}. Since $s>2\ell+2$, the same argument also gives the lower-order estimate
\[
 \sum_{\substack{|\alpha|+|\beta|\leq s-1\\0<\gamma\leq\alpha}}
 \left\|w_m(\partial_x^\gamma E)
 \partial_x^{\alpha-\gamma}\partial_v^{\beta+1}g\right\|_2
 \leq C\|\nabla E\|_{H^{s-2}_x}\|g\|_{H^{s-1}_m}.
\]
Indeed, in the high-derivative case the required number of derivatives of $g$ is at most $(s-1)-|\gamma|+1+\ell\leq2\ell\leq s-1$.
We use this estimate for the commutators in the difference equations.
For the source $E[h]\cdot\nabla_v g$, the term with no derivative on $E[h]$ is estimated separately by $\|E[h]\|_\infty\|g\|_{H^s_m}$.

Consider the linear transport equation
\[
 \partial_tg+v\cdot\nabla_xg+E(t,x)\cdot\nabla_vg=0.
\]
Apply $\partial_x^\alpha\partial_v^\beta$ with
$|\alpha|+|\beta|\leq s$. The commutator with free transport is
\[
 [\partial_x^\alpha\partial_v^\beta,v\cdot\nabla_x]g
 =\sum_{i=1}^d\beta_i
 \partial_x^{\alpha+e_i}\partial_v^{\beta-e_i}g.
\]
These terms have the same total order as the original derivative.
The commutator with the field is
\[
 [\partial_x^\alpha\partial_v^\beta,E\cdot\nabla_v]g
 =\sum_{0<\gamma\leq\alpha}
 \binom{\alpha}{\gamma}(\partial_x^\gamma E)\cdot
 \partial_x^{\alpha-\gamma}\partial_v^{\beta+1}g.
\]
Estimate this sum by \eqref{eq:A-product}.

The transport vector field $(v,E)$ is divergence-free in $(x,v)$.
When we multiply a differentiated equation by $w_m^2\partial_x^\alpha\partial_v^\beta g$ and integrate, the only term produced by the weight is bounded by
\[
 \|E\|_\infty\|w_m\partial_x^\alpha
 \partial_v^\beta g\|_2^2,
 \qquad |\nabla_v w_m|\leq C_mw_m.
\]
Summing over $\alpha,\beta$ gives
\begin{equation}\label{eq:A-linear-energy}
 \frac{d}{dt}\|g(t)\|_{H^s_m}
 \leq C\bigl(1+\|E(t)\|_\infty
 +\|\nabla E(t)\|_{H^{s-1}_x}\bigr)
 \|g(t)\|_{H^s_m}.
\end{equation}
For the nonlinear field $E[g]$, mass is conserved. Substituting the field bounds into \eqref{eq:A-linear-energy} gives
\[
 \frac{d}{dt}\|g(t)\|_{H^s_m}
 \leq C\bigl(1+M+\|g(t)\|_{H^s_m}\bigr)
 \|g(t)\|_{H^s_m},
 \qquad M=\|g(0)\|_1.
\]
This gives a common short-time bound whenever the initial $H^s_m$ norm and the mass are uniformly bounded.

\subsection{\texorpdfstring{Proof of Lemma~\ref{lem:stability}}{Proof of Lemma 2.1}}

We use a transport iteration. Set $g^{(0)}=g_{\rm in}$ and, for $n\geq0$, solve
\[
 \partial_tg^{(n+1)}+v\cdot\nabla_xg^{(n+1)}
 +E[g^{(n)}]\cdot\nabla_vg^{(n+1)}=0,
 \qquad g^{(n+1)}(0)=g_{\rm in}.
\]
The estimates in A.1 bound each field and its first derivatives.
The corresponding linear transport equation therefore has a well-defined characteristic flow on a common short time interval. It preserves nonnegativity and the $L^1$ norm. More precisely, \eqref{eq:A-linear-energy} gives
\[
 \frac{d}{dt}\|g^{(n+1)}(t)\|_{H^s_m}
 \leq C\bigl(1+M+\|g^{(n)}(t)\|_{H^s_m}\bigr)
 \|g^{(n+1)}(t)\|_{H^s_m},
 \qquad M=\|g_{\rm in}\|_1.
\]
Choose a bound larger than twice the initial norm. If the $n$-th iterate stays below this bound, the displayed differential inequality and Gronwall's inequality show that the next iterate also stays below it, provided $T$ is sufficiently small. This closes the induction and gives a common $H^s_m$ bound, independent of $n$.

We also need a uniform $W^{1,1}$ bound. Differentiating the linear equation in $x_i$ and $v_i$ gives
\[
 \begin{aligned}
 (\partial_t+v\cdot\nabla_x+E[g^{(n)}]\cdot\nabla_v)
 \partial_{x_i}g^{(n+1)}
 &=-\partial_{x_i}E[g^{(n)}]\cdot\nabla_vg^{(n+1)},\\
 (\partial_t+v\cdot\nabla_x+E[g^{(n)}]\cdot\nabla_v)
 \partial_{v_i}g^{(n+1)}
 &=-\partial_{x_i}g^{(n+1)}.
 \end{aligned}
\]
The transport field is divergence-free. Taking $L^1$ norms and summing these inequalities yields
\[
 \frac{d}{dt}\|\nabla_{x,v}g^{(n+1)}\|_1
 \leq C\bigl(1+\|\nabla E[g^{(n)}]\|_\infty\bigr)
 \|\nabla_{x,v}g^{(n+1)}\|_1.
\]
The $H^s_m$ bound and Sobolev embedding control $\|\nabla E[g^{(n)}]\|_\infty$. Gronwall's inequality therefore gives a uniform $W^{1,1}$ bound for the iterates.

We next prove convergence. Let $h^{(n+1)}=g^{(n+1)}-g^{(n)}$. For $n\geq1$, subtraction of two successive equations gives
\[
 (\partial_t+v\cdot\nabla_x+E[g^{(n)}]\cdot\nabla_v)h^{(n+1)}
 =-E[h^{(n)}]\cdot\nabla_vg^{(n)},
 \qquad h^{(n+1)}(0)=0.
\]
The field estimates imply
\[
 \|E[h]\|_\infty\leq C\bigl(\|h\|_1+\|h\|_{H^{s-1}_m}\bigr),
 \qquad
 \|\nabla E[h]\|_{H^{s-2}_x}
 \leq C\|h\|_{H^{s-1}_m}.
\]
For the source term, the part with no derivative on $E[h^{(n)}]$ is bounded by $\|E[h^{(n)}]\|_\infty\|g^{(n)}\|_{H^s_m}$.
For terms with at least one derivative on the field, use the order $s-1$ version of \eqref{eq:A-product}, with $\nabla E[h^{(n)}]$ in $H^{s-2}_x$ and $g^{(n)}$ bounded in
$H^s_m$. Together with the field estimates, this bounds the source in $H^{s-1}_m$ by $C X_n(t)$. The commutator with the transport field $E[g^{(n)}]$ is bounded by $C X_{n+1}(t)$, using the same lower-order product estimate and the uniform $H^s_m$ bound. The $L^1$ source is bounded by $\|E[h^{(n)}]\|_\infty\|\nabla_v g^{(n)}\|_1$, and the latter factor is uniformly bounded. We therefore obtain
\[
 \frac{d}{dt}X_{n+1}(t)
 \leq C X_{n+1}(t)+C X_n(t),
 \qquad
 X_n(t)=\|h^{(n)}(t)\|_1+\|h^{(n)}(t)\|_{H^{s-1}_m}.
\]
Since $X_{n+1}(0)=0$, Gronwall's inequality gives
\[
 \sup_{0\leq t\leq T}X_{n+1}(t)
 \leq C T e^{CT}\sup_{0\leq t\leq T}X_n(t).
\]
After shortening $T$, the factor on the right is less than one. The first difference is bounded by the uniform estimates. The later differences form a geometric series, so the iterates converge in $C([0,T];L^1\cap H^{s-1}_m)$. The field estimate shows that their fields converge uniformly. Since $m>1$, convergence in $H^{s-1}_m$ also gives $v\cdot\nabla_xg^{(n+1)}\to v\cdot\nabla_xg$ in $L^2_{x,v}$. Moreover, $\nabla_vg^{(n+1)}\to\nabla_vg$ strongly in $L^2_{x,v}$. Together with the uniform convergence of the fields, this shows that $E[g^{(n)}]\cdot\nabla_vg^{(n+1)}$ converges in $L^2_{x,v}$. We may therefore pass to the limit in the transport equation. The limit solves \eqref{eq:VP} and retains the $H^s_m$ bound by weak lower semicontinuity. Nonnegativity passes to the limit by strong $L^2$ convergence, and strong $L^1$ convergence preserves the mass. Since $s>d+4$, Sobolev embedding gives local $C^1$ regularity in $(x,v)$, and the equation gives a continuous time derivative. Thus the solution is classical.

It also retains the $W^{1,1}$ bound. Fix $t\in[0,T]$. Convergence in $H^{s-1}_m$ gives strong $L^2$ convergence of the function and all its first derivatives. Choose a subsequence along which these finitely many quantities converge almost everywhere. Fatou's lemma then gives
\[
 \|g(t)\|_{W^{1,1}}
 \leq\liminf_{j\to\infty}
 \|g^{(n_j)}(t)\|_{W^{1,1}}.
\]
The bound on the right is independent of $t$. Thus the limit has the required uniform $W^{1,1}$ bound.

\par\medskip
We next compare two solutions and prove the stability estimate.

Let $g$ and $\widetilde g$ be two solutions with the uniform bounds above, and put $h=g-\widetilde g$. Their difference satisfies
\[
 \partial_th+v\cdot\nabla_xh+E[g]\cdot\nabla_vh
 =-E[h]\cdot\nabla_v\widetilde g.
\]
For the source term, the product estimate gives
\[
 \|E[h]\cdot\nabla_v\widetilde g\|_{H^{s-1}_m}
 \leq C_T\bigl(\|E[h]\|_\infty
 +\|\nabla E[h]\|_{H^{s-2}_x}\bigr)
 \leq C_T\bigl(\|h\|_{H^{s-1}_m}+\|h\|_1\bigr).
\]
The commutator with the transport field $E[g]$ is bounded by $C_T\|h\|_{H^{s-1}_m}$, by the order $s-1$ product estimate and the uniform $H^s_m$ bound for $g$. Hence
\[
 \frac{d}{dt}\|h\|_{H^{s-1}_m}
 \leq C_T\bigl(\|h\|_{H^{s-1}_m}+\|h\|_1\bigr).
\]
For the $L^1$ estimate, the divergence-free transport field gives
\[
 \frac{d}{dt}\|h\|_1
 \leq \|E[h]\|_\infty\|\nabla_v\widetilde g\|_1.
\]
The $W^{1,1}$ estimate proved above bounds the last factor uniformly.
Using the field estimate once more, we find
\[
 \frac{d}{dt}\|h\|_1
 \leq C_T\bigl(\|h\|_1+\|h\|_{H^{s-1}_m}\bigr).
\]
Adding the two inequalities and applying Gronwall yields
\begin{equation}\label{eq:A-stability}
 \sup_{0\leq t\leq T}
 \left(\|g(t)-\widetilde g(t)\|_1
 +\|g(t)-\widetilde g(t)\|_{H^{s-1}_m}\right)
 \leq C_T\left(\|g(0)-\widetilde g(0)\|_1
 +\|g(0)-\widetilde g(0)\|_{H^{s-1}_m}\right).
\end{equation}
In particular, this estimate gives uniqueness in the class constructed above.

Apply \eqref{eq:A-stability} to the solutions with initial data $g_{\rm in}^\varepsilon$ and $g_{\rm in}^0$. The assumed convergence
of the initial data gives
\[
 \sup_{0\leq t\leq T}
 \|g^\varepsilon(t)-g^0(t)\|_{H^{s-1}_m}\longrightarrow0.
\]
By \eqref{eq:A-density},
\[
 \sup_{0\leq t\leq T}
 \|\rho^\varepsilon(t)-\rho^0(t)\|_{H^{s-1}_x}
 \longrightarrow0.
\]
Since $s-1>d/2+2$, Sobolev embedding gives
\[
 \rho^\varepsilon\longrightarrow\rho^0
 \quad\text{in }C([0,T];C^2_{\rm loc}(\R^d)),
\]
which is the convergence claimed in Lemma~\ref{lem:stability}.

Finally, we check the time regularity used in Section~2.1. We first justify the moment equations. Choose $\chi\in C_c^\infty(\R^d)$ such that $\chi=1$ on $B_1$ and 
$\chi=0$ outside $B_2$, and set $\chi_R(v)=\chi(v/R)$. Test the weak Vlasov equation with $\psi(t,x)\chi_R(v)$ and $\psi(t,x)v_i\chi_R(v)$, where
$\psi\in C_c^\infty((0,T)\times\R^d)$. The weighted $L^2$ bound gives, for every compact $K\subset\R^d$,
\[
 \int_K\int_{|v|\geq R}\la v\ra^2|g(t,x,v)|\,dv\,dx
 \leq C_K
 \left(\int_{|v|\geq R}\la v\ra^{4-2m}\,dv\right)^{1/2}
 \|\la v\ra^m g(t)\|_{L^2_{x,v}},
\]
which tends to zero as $R\to\infty$ because $m>d/2+2$. The same estimate holds after integration in time. Since
$|\nabla_v\chi_R|\leq C/R$ and
$|v_i\nabla_v\chi_R|\leq C$ on the annulus
$R\leq|v|\leq2R$, all cutoff error terms tend to zero. Passing to the
limit gives, in distributions,
\[
 \partial_t\rho+\nabla_x\cdot j=0,
 \qquad
 \partial_tj+\nabla_x\cdot\Pi=\rho E.
\]
Thus both moment equations are justified by the cutoff limit.

For $|\alpha|\leq s-1$, weighted Cauchy--Schwarz gives
\[
 \|\partial_x^\alpha j(t)\|_2
 +\|\partial_x^\alpha\Pi(t)\|_2
 \leq C_m
 \|\la v\ra^m\partial_x^\alpha g(t)\|_{L^2_{x,v}}.
\]
The integral is finite because $m>d/2+2$. Since $g\in C([0,T];H^{s-1}_m)$, we obtain $j,\Pi\in C([0,T];H^{s-1}_x)$.
The moment equations therefore give
\[
 \partial_t\rho=-\nabla_x\cdot j,
 \qquad
 \partial_t^2\rho
 =\nabla_x\cdot\nabla_x\cdot\Pi-\nabla_x\cdot(\rho E).
\]

We have $\rho\in C([0,T];H^{s-1}_x)$, $E\in C([0,T];L^\infty_x)$, and $\nabla E\in C([0,T];H^{s-2}_x)$. Put $k=s-2$. Since $k-1=s-3>d/2$, the space $H^{k-1}(\R^d)$ is an algebra. By the Leibniz rule,
\[
 \|\rho E\|_{H^k_x}
 \leq C\left(
 \|E\|_\infty\|\rho\|_{H^k_x}
 +\|\rho\|_{H^{k-1}_x}\|\nabla E\|_{H^{k-1}_x}
 \right).
\]
The first term contains the part with no derivative on $E$. All other terms are derivatives of order at most $k-1$ of $\rho\nabla E$ and are controlled by the algebra property. Thus $\rho E\in C([0,T];H^{s-2}_x)$, and
\[
 \partial_t\rho\in C([0,T];H^{s-2}_x),
 \qquad
 \partial_t^2\rho\in C([0,T];H^{s-3}_x).
\]
Since $d\geq2$ and $s>d+4$, we have $s-3>d/2+2$. Sobolev embedding gives
\[
 \rho,\ \partial_t\rho,\ \partial_t^2\rho
 \in C([0,T];C^2_{\rm loc}(\R^d)).
\]
In particular, $\rho\in C^2([0,T];C^2_{\rm loc}(\R^d))$.

The initial density is strictly positive. The characteristic formula for the smooth solution preserves positivity, so $\rho$ stays positive near the point used in Section~2.1. Hence $u=\log\rho$ has the same time and space regularity there. This justifies the Taylor expansion in \eqref{eq:loss-flat} and completes the proof of Lemma~\ref{lem:stability}.

\medskip
\noindent\textbf{Data/Code Availability.} The manuscript has no associated data or code.

\medskip
\noindent\textbf{Conflict of interest.} The author declares that he has no conflict of interest.

\end{document}